% To: Mijrady biyrujalyim <shlhetal@sundial.ma.huji.ac.il>
% Subject: 643 final version
% Date: Tue, 11 Jan 2000 14:31:31 +0200 (IST)
% From: Saharon Shelah <shelah@math.huji.ac.il>
% Mime-Version: 1.0
% X-sliced-and-diced-by: 'savemail' 0.3, Feb 1999

%TeX File
%
%\magnification=\magstep1

%\baselineskip24pt
\input amssym.tex
\input amssym.def
\overfullrule0pt

\centerline{
\bf CARDINAL INVARIANTS ${\frak b}_\kappa$ AND ${\frak t}_\kappa$
}

\bigskip

\centerline{$\hbox{\bf Saharon Shelah}^{\dag}$}

\centerline{ Department of Mathematics}

\centerline{ Hebrew University of Jerusalem, Jerusalem, Israel}

\centerline{ Rutgers University, New Brunswick, NJ, USA}

\centerline{ University of Wisconsin, Madison, WI, USA}

\bigskip

\centerline{$\hbox{\bf Zoran Spasojevi\' c}^{\ddag}$}

\centerline{ Department  of Mathematics}

\centerline{ Massachusetts Institute of Technology}

\centerline{ Cambridge, MA 02139, USA}

\bigskip

\bigskip

\noindent ABSTRACT:\ \  This paper studies cardinal invariants 
${\frak b}_\kappa$ and ${\frak t}_\kappa$, the natural generalizations
of the invariants ${\frak b}$ and ${\frak t}$ to a regular cardinal
$\kappa$.

\bigskip

\noindent {\bf \S1. Introduction}

\bigskip

Cardinal invariants ${\frak b}$ and ${\frak t}$ were introduced by
Rothberger [3, 4]. They are cardinals between $\omega_1$ and
$2^\omega$ and have been extensively studied over the years. 
The survey paper [5] contains much information about these two
invariants as well as many other cardinal invariants of the continuum.

The goal of this paper is to study the natural generalizations of
${\frak b}$ and ${\frak t}$ to higher regular cardinals, namely
${\frak b}_\kappa$ and ${\frak t}_\kappa$  respectively, where
$\kappa$ is a regular cardinal. The results
presented here are that the relationship ${\frak t} \leq {\frak b}$
(shown by Rothberger [4]) also holds for ${\frak b}_\kappa$ and 
${\frak t}_\kappa$  and that, under certain cardinal arithmetic
assumption, if $\kappa \leq \mu <{\frak t}_\kappa$ then 
$2^\kappa = 2^\mu$. These results are then used as constraints
in the forcing construction of model in which ${\frak b}_\kappa$
and ${\frak t}_\kappa$ can take on essentially any preassigned regular value.
%In the last section, a class forcing is used to simultaneously control
%the values of ${\frak b}_\kappa$ and ${\frak t}_\kappa$ for all
%regular cardinals $\kappa$.

\bigskip

\noindent {\bf \S2 Conventions and elementary facts}

\bigskip

For cardinals $\lambda$ and $\kappa$  let
$[\kappa]^\lambda = \{ X \subseteq \kappa : \mid X \mid = \lambda \}$
and $\kappa^\lambda = \{ f : f: \lambda \to \kappa \}$.
The symbol $\kappa^\lambda$ is also used to denote the cardinality
of the set $\{ f : f : \lambda \to \kappa \}$ but the intended meaning 
will be clear from the context. For $A,B \in [\kappa]^\kappa$
let $A \subseteq^* B$ iff $\mid A \setminus B \mid < \kappa$ and
$A \subset^* B$ iff $\mid A \setminus B\mid < \kappa \ \wedge
		     \mid B \setminus A\mid = \kappa$.
For $f,g \in \kappa^\kappa$ let $f <^* g$ iff 
$\exists \beta < \kappa \forall \alpha > \beta 
	      ( f(\alpha) < g(\alpha) )$.
Then ${\cal B} \subseteq \kappa^\kappa$ is unbounded in 
     $(\kappa^\kappa,<^*)$ if
     $\forall f \in \kappa^\kappa \exists g \in {\cal B}
      (g \not <^* f)$.
Let

${\frak b}_\kappa = {\rm min} \{ \mid {\cal B} \mid : \hbox{
 ${\cal B} \subset \kappa^\kappa$ and ${\cal B}$ is $<^*$-unbounded
  in $\kappa^\kappa$ }\}$,

${\frak t}_\kappa = {\rm min} \{ \mid {\cal T} \mid : \hbox{
 ${\cal T} \subseteq [\kappa]^\kappa$, $\mid {\cal T} \mid \geq \kappa$,
 ${\cal T}$ is well ordered by $\subset^*$,}$

\hskip 2.4cm $\hbox{
$\forall C\in [\kappa]^\kappa ( \mid \kappa \setminus C \mid = \kappa
     \to \exists A \in {\cal T} ( \mid A \setminus C \mid = \kappa ))$ }\}$.

In this notation ${\frak b} = {\frak b}_\omega$ and 
${\frak t} = {\frak t}_\omega$. An equivalent formulation of 
${\frak t}_\kappa$ is obtained if $\subseteq^*$ is used instead of
$\subset^*$. Standard arguments show that
$ \kappa^+ \leq {\frak b}_\kappa , {\frak t}_\kappa \leq 2^\kappa$
and that in the definition of ${\frak b}_\kappa$, ${\cal B}$ may
be assumed to be well ordered by $<^*$ and consisting only
of strictly increasing functions. Thus, both ${\frak b}_\kappa$
and ${\frak t}_\kappa$ are regular cardinals.

\bigskip

\hrule

\bigskip

\noindent {\dag} Partially supported by the NSF and Israel Science Foundation,
founded by the Israel Academy of Sciences. Pub. 643

\medskip

\noindent {\ddag} Supported by the NSF Grant No. DMS-9627744

\vfill\eject

\noindent {\bf Lemma 1.}\ \  ${\frak t}_\kappa \leq {\frak b}_\kappa$

\bigskip

\noindent Proof:\ \ The case $\kappa = \omega$ was established
in [4]. So assume $\kappa > \omega$ and by
way of contradiction assume ${\frak b}_\kappa < {\frak t}_\kappa$.
Let $\{ f_\alpha : \alpha < {\frak b}_\kappa \} \subseteq
\kappa^\kappa$ be $<^*$-unbounded in $(\kappa^\kappa, <^*)$ 
and such that $\alpha < \beta \to f_\alpha <^* f_\beta$.
For each $\alpha < {\frak b}_\kappa$ let
$C_\alpha = \{ \xi < \kappa : \forall \eta < \xi
			     ( f_\alpha (\eta) < \xi) \}$.
Then each $C_\alpha$ is closed unbounded in $\kappa$ and
$\alpha \leq \beta \to C_\beta \subseteq^* C_\alpha$.
Since ${\frak b}_\kappa < {\frak t}_\kappa$,
$\exists A \in [\kappa]^\kappa \forall \alpha < {\frak b}_\kappa
 ( A \subset^* C_\alpha)$. Let $f : \kappa \to A$ be such that
$\forall \xi < \kappa (\xi < f(\xi))$. Fix $\alpha < {\frak b}_\kappa$
and let $i_\alpha$ be such that 
$A \setminus i_\alpha \subseteq C_\alpha \setminus i_\alpha$.
Let $\xi \in \kappa \setminus i_\alpha$ and
$\eta = {\rm min} (C_\alpha \setminus (\xi + 1))$ and note that
$f_\alpha (\xi) < \eta$. However, 
$A \setminus \xi \subseteq C_\alpha \setminus \xi$, and $\xi < f (\xi )$,
so $\eta \leq f (\xi)$. In other words,
$\forall \xi \in \kappa \setminus i_\alpha (f_\alpha (\xi) < f(\xi))$
so that $f$ is a bound for $\{ f_\alpha : \alpha < {\frak b}_\kappa \}$.
This is a contradiction and the lemma is proved. \qquad $\square$

\bigskip

The cardinal ${\frak b}_\kappa$ was studied in [1] where it was shown
that the value of ${\frak b}_\kappa$ does not have any influence 
on the value of $2^\mu$ for $\kappa \leq \mu < {\frak b}_\kappa$
even if $GCH$ is assumed to hold below $\kappa$.
However, the same does not hold for ${\frak t}_\kappa$ as it is shown
in the next section.

\bigskip

\noindent{\bf \S3 Combinatorics}

\bigskip

The goal of this section is to show that if $\kappa^{<\kappa} = \kappa$
then $ 2^\mu = 2^\kappa$ for any regular $\mu$ with
$\kappa \leq \mu < {\frak t}_\kappa$. The idea behind the proof
is essentially the same as that of the proof of
$\omega \leq \mu < {\frak t} \to 2^\mu = 2^\omega$, namely
to use $\mu < {\frak t}_\kappa$ to construct a binary tree
in $({\cal P}(\kappa) , \subset^*)$ of hight $\mu$. However,
unlike in the case $\kappa = \omega$, when $\kappa$ is uncountable
a difficulty arises in the construction at limit stages of cofinality
less than $\kappa$. The difficulty comes from the fact that the intersection
of a $\subset^*$-decreasing sequence in $[\kappa]^\kappa$ of limit
length less than $\kappa$ may be empty. To deal with this difficulty,
a notion of a closed subset of $\kappa$ with respect to a certain
parameter is introduced next.

Let $D$  be a filter on a regular cardinal $\kappa$ and
$A,B \in [\kappa]^\kappa$. Then $ A=_D 0$ if $\kappa \setminus A \in D$
and $A =_D B$ if $(A \setminus B)\cup (B \setminus A) =_D 0$. Let $D^+$
be the collection of all sets $A \subseteq \kappa$ such that 
$\kappa \setminus A \not\in D$. If $A,B \in D^+$ then
$A \subset^*_D B$ if $A \setminus B =_D 0$ and $B \setminus A \in D^+$.
For each $i < \kappa$ let $A^\kappa_i \subseteq \kappa \setminus (i+1)$
be such that $3(i+1) \in A^\kappa_i$, 
$\mid A^\kappa_i \mid = \kappa$, 
$\bigcup_{i<\kappa} A^\kappa_i = \kappa \setminus \{ 0 \}$, and
if $i < j < \kappa$ then $A^\kappa_i \cap A^\kappa_j = \emptyset$.
Let $E$ be the set of all limit ordinals $\delta < \kappa$
such that $\delta$ is a cardinal or it is a multiple of 
$\mid \delta \mid^\omega$ (ordinal exponentiation) and for each 
$\alpha < \delta$, $A^\kappa_\alpha \cap \delta $ is unbounded 
in $\delta$. Let $D^\delta$ be the collection of all subsets $X$ of 
$\delta$ such that for some $\alpha < \delta$,
$([\alpha, \delta) \setminus \bigcup_{i < \alpha} A^\kappa_i)
\subseteq X$ (we are interested in the cases when $\delta \in E$). 
And let $D(\kappa)$ be the collection of all subsets
$A$ of $\kappa$ such that for some $\alpha < \kappa$,
$([\alpha, \kappa) \setminus \bigcup_{i < \alpha} A^\kappa_i)
\subseteq A$. 
Let $A^\kappa = \bigcup \{ A^\kappa_\delta : \delta \in E\}$ and
for $\zeta \in A^\kappa$ choose $D_\zeta$ and $\delta_\zeta$
such that $\zeta \in A^\kappa_{\delta_\zeta}$, $D_\zeta$ is
a filter on $\delta_\zeta$ generated by less than $\kappa$ sets
which extends
$D^{\delta_\zeta}$, and for any $\delta \in E$ and any
filter $D$ on $\delta $ which extends $D^\delta$ and which is
generated by less than $\kappa$, 
$\mid \{ \zeta \in A^\kappa_\delta : D = D_\delta \}\mid = \kappa$.
Let $\bar D = \langle D_\zeta : \zeta \in A^\kappa \rangle$.

\bigskip

\noindent {\bf Definition 2.}\ \ A subset $A$ of $\kappa$ is
$( E,\bar D )$-closed if for every $\delta \in E$
and $\zeta \in A^\kappa_\delta$ then $\zeta \in A$ whenever
$A\cap \delta \in D_\zeta$. Let $cl^0(A) = A$,
$cl^1(A) = A \cup \{ \zeta \in A^\kappa : A \cap \delta_\zeta \in D_\zeta\}$
and $cl^\alpha (A) = A \cup \{ cl^1(cl^\beta (A)) : \beta < \alpha \}$.
Let $(E, \bar D )$-closure of $A$, $cl(A)$, be the set
$cl^\alpha (A)$ for every $\alpha$ large enough.

\bigskip

The above definition formulates the notion of a closed set with respect
to a parameter $(E,\bar D)$. The following sequence of observations gives
some elementary properties of the closure operations which will be
needed in the sequel.

\bigskip

\noindent {\bf Observation 3.}\ \  Let  $A \subseteq B \subseteq \kappa$.
Then

$(a)$\ \ $\kappa$ is an $(E, \bar D)$-closed subset of $\kappa$,

$(b)$\ \ $\forall \alpha 
      ( cl^\alpha (A) \subseteq cl^\alpha (B) \subseteq \kappa )$,

$(c)$\ \ $cl(A)$ is the minimal $(E,\bar D)$-closed set which contains $A$,

$(d)$\ \ $A =_{D(\kappa)} 0$ iff $cl(A) =_{D(\kappa)} 0$,

$(e)$\ \ if for some $\zeta \in A^\kappa_\delta$, 
	 $\zeta \in cl^1 (\delta \cap A)$ then
	 $\mid cl(A) \cap A^\kappa_\delta \mid = \kappa$,

\bigskip

The next lemma is used in the construction of the successor levels of 
the binary tree.

\bigskip

\noindent {\bf Lemma 4.}\ \ If $A$ is $(E, \bar D)$-closed and from
$D^+(\kappa)$ then there are two disjoint $(E,\bar D)$-closed subsets
of $A$ in $D^+(\kappa)$.

\bigskip

\noindent Proof:\ \  By induction on $i < \kappa$ choose distinct 
ordinals $\alpha_i , \beta_i \in 
	   A \setminus cl( {\rm sup}\{ \alpha_j + 1, \beta_j + 1 : j<i\})$.
Then $cl (\{ \alpha_i : i < \kappa \})$ and
 $cl ( \{ \beta_i : i < \kappa \})$ are as desired. \qquad $\square$

\bigskip

The following lemma is used in the construction of the limit levels
of the tree whose cofinality is less than $\kappa$.

\bigskip

\noindent {\bf Lemma 5} \ \  Let $\tau < \kappa$ be a regular cardinal.
Let $A_i \in D^+(\kappa)$ $(i<\tau)$ be $(E, \bar D)$-closed such
that $i < j < \kappa \to A_j \subset^*_{D(\kappa)} A_i$.
Then $\bigcap_{i < \tau} A_i \in D^+(\kappa)$ and is 
$(E, \bar D)$-closed.

\bigskip

\noindent Proof:\ \ For each $i < \tau$ let $C_i \subseteq E$ be closed
unbounded in $\kappa$ such that 
$\forall \delta \in C_i (A_i \cap \delta \not =_{D^\delta} 0)$.
Then $C = \bigcap_{i<\tau} C_i$ is also closed unbounded in $\kappa$.
For each $\delta \in C$ let $D^*_\delta$ be the filter on $\delta$
generated by $D^\delta\cup \{ A_I\cap \delta : i<\tau \}$.
Clearly $\delta \in C \wedge i<\tau \to A_i \cap \delta \in D^*_\delta$,
hence $B_\delta = \{ \zeta : \zeta\in A^\kappa_\delta \wedge
		     D_\zeta = D^*_\delta \}$ is an unbounded
subset of $A^\kappa_\delta$ and is a subset of $A_i$ for each
$i < \tau$, since each $A_i$ is $(E,\bar D)$-closed.
Then $\bigcup_{\delta \in C} B_\delta \in D^+(\kappa)$
witnesses that $\bigcap_{i<\tau} A_i$ is as required. \qquad $\square$

\bigskip

And the final lemma of this section will aid in the construction of the limit
levels of the tree of cofinality greater than or equal to $\kappa$.

\bigskip

\noindent {\bf Lemma 6.}\ \  Let $\tau$ be a regular cardinal with
$\kappa \leq \tau < {\frak t}_\kappa$ and
$\langle A_i : i < \tau \rangle \subseteq D^+(\kappa)$ such
that each $A_i$ is $(E,\bar D)$-closed and 
$i < j < \tau \to A_j \subset^*_{D(\kappa)} A_i \wedge
 A_i \not \subset^*_{D(\kappa)} A_j$.
Then there is an $(E,\bar D)$-closed $B \in D^*(\kappa)$
such that $\forall i < \tau ( B \subset^*_{D(\kappa)} A_i)$.

\bigskip

\noindent Proof:\ \  For each $i < j < \tau$ let $C_{ij} \subseteq E$
be closed unbounded in $\kappa$ such that $A_i \cap \delta$,
$A_j \cap \delta$, $(A_i \setminus A_j) \cap \delta \in (D^\delta)^+$
for each $\delta \in C_{ij}$. Let $g_{ij}$ be the function
enumerating $C_{ij}$ in the increasing order.
Since ${\frak t}_\kappa \leq {\frak b}_\kappa$ let
$g : \kappa \to E$ be a strictly increasing function such that
$g_{ij} <^* g$ for each $i,j < \tau$. Let $C$ be the collection
of all limit points of ${\rm ran}(g)$. Then $C \subseteq E$ is closed
unbounded in $\kappa$ and $C \subset^* C_{ij}$ for each $i,j < \tau$.
Now, for $i,j < \tau$, define $f_{ij} : C \to \kappa$ by:
$f_{ij} ( \delta )$ is the least $\zeta \in A^\kappa_\delta$ such that 
$(A_i \setminus A_j) \cap \delta \in D_\zeta$ if such $\zeta$
exists (and it does exist whenever $\delta \in C_{ij}$) and zero
otherwise. And again, since ${\frak t}_\kappa \leq {\frak b}_\kappa$,
let $f : C \to \kappa$ be such that for each $i < j < \tau$,
$f_{ij} <^* f$. Now let
$X = \bigcup \{ f(\delta) \cap A^\kappa_\delta : \delta \in C \}$
and $A'_i = A_i \cap X$ for $i < \tau$ and note that $A'_i$ is
an unbounded subset of $\kappa$ and
$i < j < \tau \to A'_j \subset^* A'_i \wedge
    A'_i \not \subset^* A'_j$.
Then, by the definition of ${\frak t}_\kappa$, since 
$\tau < {\frak t}_\kappa$, there is a $B^* \subseteq\kappa$,
an unbounded subset of $X$, such that
$i < \tau \to B^* \subset^* A'_i$. In addition, $B^* \in D^+(\kappa)$
and by the choice of $X$ and
$B^*\subseteq A'_{i+1}\subseteq A_{i+1} \subset^*_{D(\kappa)} A_i$
for each $i < \tau$.
Hence, by Observation 3(b), 
$B = cl (B^*)\subseteq^*_{D(\kappa)} cl (A_{i+1}) = A_{i+1}
\subset^*_{D(\kappa)} A_i$ is as desired
and the lemma is proved. \qquad $\square$

\bigskip

At this point enough preliminary work is completed for the proof
of the main result of this section.

\bigskip

\noindent {\bf Theorem 7.}\ \ Let $\kappa$ and $\mu$ be a regular cardinals
such that $\kappa^{<\kappa} = \kappa$ and 
$\kappa \leq \mu < {\frak t}_\kappa$. Then $2^\mu = 2^\kappa$.

\bigskip Proof:\ \  By induction on $\zeta \leq \mu$, for every sequence
$\eta$ of zeros and ones of length $\zeta$, choose a set $A_\eta$
such that

$(a)$\ \ $A_\eta \in D^+(\kappa)$,

$(b)$\ \ $A_\eta$ is $(E, \bar D)$-closed,

$(c)$\ \ if $\rho$ is an initial segment of $\eta$ then
    $A_\eta \setminus A_\rho =_{D(\kappa)} 0$ and
    $A_\rho \setminus A_\eta \in D^+(\kappa) $,

$(d)$\ \ $A_{\eta^\frown\langle 0 \rangle} \cap
      A_{\eta^\frown\langle 1\rangle} = \emptyset$.

Let $A_{\langle \emptyset \rangle} = \kappa$. For the successor step
suppose $\eta$ is a sequence of zeros and ones and $A_\eta$
has been constructed. By Lemma 4, let $B$ and $C$ be two disjoint 
$(E, \bar D)$-closed subsets of $A_\eta$ from $D^+(\kappa)$.
Let $A_{\eta^\frown \langle 0 \rangle} = B$
and $A_{\eta^\frown \langle 1 \rangle} =C$. This takes care of the successor
stages of the construction.

Now suppose $\eta$ is a sequence of zeros and ones such that
${\rm dom} (\eta) = \lambda \leq \mu$ is a limit ordinal with
${\rm cf} (\lambda) = \tau < \kappa$, and for each $\alpha < \lambda$,
$A_{\eta \upharpoonright \alpha}$ has been constructed.
Let $\langle \alpha_i : i < \tau \rangle$ be an increasing sequence
of ordinals with limit $\lambda$. By Lemma 5,
$\bigcap_{i < \tau} A_{\eta \upharpoonright \alpha_i} \in D^+(\kappa)$
and is $(E, \bar D)$-closed. Let
$A_\eta = \bigcap_{i<\tau} A_{\eta\upharpoonright \alpha_i}$.
This takes care of limit stages of cofinality less than $\kappa$.

Finally suppose $\eta$ is a sequence of zeros and ones such that
${\rm dom} (\eta) = \lambda \leq \mu$ is a limit ordinal with
$\kappa \leq {\rm cf}(\lambda) = \tau$ and for each $\alpha < \lambda$,
$A_{\eta\upharpoonright\alpha}$ has been defined.
Let  $\langle \alpha_i : i < \tau \rangle$ be an increasing sequence
of ordinals with limit $\lambda$. By Lemma 6, there is an
$(E,\bar D)$-closed $B \in D^+(\kappa)$ such that
$\forall i < \tau 
 (B \subset^*_{D(\kappa)} A_{\eta\upharpoonright\alpha_i})$.
Let $A_\eta = B$. Now $\{ A_\eta : \eta \in 2^\mu \}$
is a family of $2^\mu$ distinct subsets of $\kappa$ so that
$2^\mu \leq 2^\kappa$ and the proof is finishes. \qquad $\square$

\bigskip

\noindent {\bf \S4 Forcing}

\bigskip

Let $\kappa$ be a regular uncountable cardinal and let
$\lambda$, $\mu$, $\theta$ be cardinals such that
$\kappa < \lambda \leq \mu \leq \theta$ with $\lambda$, $\mu$
regular and ${\rm cf}(\theta) > \kappa$.
This section deals with the construction of a model for
${\frak t}_\kappa = \lambda$, ${\frak b}_\kappa = \mu$
and $2^\kappa = \theta$. The idea behind the construction is as
follows: Start with a countable transitive model (c.t.m.) $N$
for $ZFC + GCH$. Expend $N$ to a model $M$ by using the
standard partial order for adding $\theta^+$ 
many subsets of $\lambda$ (see below).
Then 
$$M\models\hbox{ ``\ $\forall \xi < \lambda (2^\xi = \xi^+ \wedge
	       2^\lambda = \theta^+)$ \ ''}.$$
In $M$, perform an iterated forcing construction 
with $<\!\!\kappa$-supports of length
$\theta\cdot \mu$ (ordinal product) with $\kappa$-closed and 
$\kappa^+$-cc partial orders as follows:
At stages which are not of the form $\theta \cdot \xi$ ($\xi < \mu$)
towers in $({\cal P}(\kappa),\subset^*)$ of hight
$\eta$ are destroyed for $\kappa < \eta < \lambda$. At stages of the form
$\theta \cdot \xi$  a function from $\kappa$ to $\kappa$ is added
to eventually dominate all the functions from $\kappa$ to $\kappa$
constructed by that stage. The bookkeeping is arranged in such a way
that by the end of the construction  all towers of
hight $\eta$ for $\kappa < \eta < \lambda$ 
are considered so that in the final model
${\frak t}_\kappa \geq \lambda$. However, in the final model
$$\forall \xi ((\xi < \kappa \to 2^\xi = \xi^+) \wedge
   (\kappa \leq \xi < \lambda \to 2^\xi = \theta)) \wedge 2^\lambda =
   \theta^+$$
so that, by the previous section, ${\frak t}_\kappa = \lambda$.
By virtue of adding dominating functions at stages of the form
$\theta \cdot \xi$, the final model has a scale in 
$(\kappa^\kappa, <^*)$ of order type $\mu$ so that 
${\frak b}_\kappa = \mu$.

The rest of this section deals with the
details of the construction.
In showing that the final model has the desired properties
it is important to know that cardinals are not collapsed.
A standard way of proving this is to show that the final
partial order obtained by the iteration is $\kappa$-closed
and has the $\kappa^+$-cc.
And to show that the final partial order has the two properties,
the names for the partial orders
used in the iteration must be carefully selected. 
The discussion here will be analogous to the discussion
in the final section of [2] which deals with countable support
iterations. Also many proofs are omited here since they are analogous to 
the proofs of the corresponding facts in [2].

\bigskip

\noindent {\bf Definition 8.}\ \  Let ${\Bbb P}$ be a partial order
and $\pi$ a ${\Bbb P}$-name for a partial order. $\pi$ is full for
$<\!\! \kappa$-sequences iff whenever $\alpha < \kappa$, $p \in {\Bbb P}$,
$\rho_\xi \in {\rm dom}(\pi)$ ($\xi < \alpha$) and for each
$\xi < \zeta < \alpha$
$$p\Vdash\hbox{`` \ $\rho_\zeta, \rho_\xi \in \pi \wedge
     \rho_\zeta \leq \rho_\xi$ \ ''}$$
then there is a $\sigma \in {\rm dom}(\pi)$ such that
$p\Vdash \hbox{`` \  $\sigma \in \pi$ \ ''}$ and
$p\Vdash \hbox{`` \ $\sigma \leq \rho_\xi$ \ ''}$ for all $\xi < \alpha$.

\bigskip

The reason for using names which are full for $<\!\!\kappa$-sequences
is because of the following

\bigskip

\noindent{\bf Lemma 9.}\ \ Let $M$ be a c.t.m. for $ZFC$ and in $M$
let
$$\langle\langle {\Bbb P}_\xi : \xi \leq \alpha \rangle,
	 \langle \pi_\xi : \xi < \alpha \rangle\rangle$$
be a $<\!\!\kappa$-support iterated forcing construction and suppose
that for each $\xi$, the ${\Bbb P}_\xi$-name $\pi_\xi$ is full
for $<\!\!\kappa$-sequences. Then ${\Bbb P}_\alpha$ is $\kappa$-closed in $M$.

\bigskip

The next few paragraphs show how to select names for partial 
orders in the construction
so that they are full for $<\!\!\kappa$-sequences. First consider
the partial order which destroys a tower in 
$({\cal P}(\kappa),\subset^*)$. Let $\epsilon$ be a regular cardinal
with $\kappa < \epsilon < \lambda$ and
$a = \langle a_\xi : \xi < \epsilon \rangle$ a tower in
$({\cal P}(\kappa),\subset^*)$. In the following subsets
of $\kappa$ are identified with their characteristic functions.

\bigskip

\noindent {\bf Definition 10.}\ \ 
${\Bbb T}_a = \{ (s,x): \hbox{ $s$ is a function} \wedge 
{\rm dom}(s) \in \kappa \wedge {\rm ran}(s) \subseteq 2 \wedge
x \in [\epsilon]^{<\kappa} \}$

\noindent with  $(s_2,x_2) \leq (s_1,x_1)$ iff

\noindent 1)\ \ $s_1 \subseteq s_2 \wedge x_1 \subseteq x_2$,

\noindent 2)\ \ $\forall \xi \in x_1 \forall \eta \in
		 {\rm dom}(s_2) \setminus {\rm dom}(s_1)
		 (a_\xi(\eta) \leq s_2(\eta))$.

\bigskip

\noindent Then ${\Bbb T}_a$ is a partial order and it is $\kappa$-closed
and $\kappa^+$-cc (assuming $\kappa^{<\kappa} = \kappa$).
Let $G$ be ${\Bbb T}_a$-generic over $M$ and 
$b = \cup \{ s : \exists x ((s,x) \in G)\}$.
Since $G$ intersects suitably chosen dense subsets of
${\Bbb T}_a$ in $M$, then
$b \subseteq \kappa$, $\mid b \mid = \mid \kappa \setminus b \mid = \kappa$
and $\forall \xi < \epsilon (a_\xi \subseteq^* b)$ so that
$a$ ceases to be a tower in $M[G]$.

Since the $<\!\!\kappa$-support iteration is sensitive to the particular names
used for the partial orders, a suitable 
name for ${\Bbb T}_a$ is formulated next.

\bigskip

\noindent {\bf Definition 11.}\ \ Assume that ${\Bbb P} \in M$,
$( {\Bbb P}\ \hbox{ is $\kappa$-closed})^M$ and
$$ {\bf 1} \Vdash \hbox{ `` \  $\tau$ is an $\check\epsilon$-tower
			    in $({\cal P}(\kappa), \subset^*)$\ ''}.$$
A standard name for ${\Bbb T}_\tau$ is 
$\langle \sigma, \leq_\sigma, {\bf 1}_\sigma \rangle$,
where
$$\eqalign{\sigma = \{
	     \langle op(\check s, \rho), {\bf 1}_{\Bbb P} \rangle :
            &\hbox{ $s$ is a function $\wedge$ ${\rm dom}(s) \in \kappa$ 
	      $\wedge$ ${\rm ran}(s) \subseteq 2$ $\wedge$}\cr
	    &\hbox{ ${\bf 1} 
	       \Vdash \hbox{ ``\ $\rho \subseteq \tau \wedge
				     \mid \rho \mid < \kappa$\ ''}$
               $\wedge$ $\rho$ is a nice name for a subset
				   of $\tau$ }\}\cr } $$
and ${\bf 1}_\sigma = op(\check 0, \check 0)$.

\bigskip

\noindent Here $op$ is the invariant name for the ordered pair
and $\rho$ is a nice name for a subset of $\tau$ if
$$\rho = \cup\{\{ \pi\}\times A_\pi: \pi \in {\rm dom}(\tau)$$
and each $A_\pi$ is an antichain in ${\Bbb P}$.
It is irrelevant what type of name we use for $\leq_\sigma$
as long as it is forced by ${\bf 1}_{\Bbb P}$ to be the correct
partial order on ${\Bbb T}_\tau$.

In $M$, let ${\Bbb P}$, $\tau$, and $\sigma$ be as in the definition above.
Let $G$ be ${\Bbb P}$-generic over $M$ and $a = \tau_G$.
Then in $M[G]$, $\sigma_G = {\Bbb T}_a$.
In addition, $\sigma$ is full for $<\!\!\kappa$-sequences.

The dominating function partial order is considered next.
Let $F \subseteq \kappa^\kappa$. In the final construction $F$ will
be equal to $\kappa^\kappa$, but for the general discussion $F$ is
any subset of $\kappa^\kappa$.

\bigskip

\noindent {\bf Definition 12.}\ \ ${\Bbb D}_F = 
\{ (s,x) : \hbox{
   $s$ is a function $\wedge$ ${\rm dom}(s) \in \kappa $
   $\wedge$ ${\rm ran}(s) \subseteq \kappa$
   $\wedge$ $x\in [F]^{<\kappa}$} \}$

\noindent where $(s_2,x_2) \leq (s_1,x_1)$ iff

\noindent 1) \ \ $s_1 \subseteq s_2 \wedge x_1 \subseteq x_2$,

\noindent 2)\ \ $\forall f \in x_1 \forall \alpha \in
		  {\rm dom}(s_2) \setminus {\rm dom}(s_1)
		   (f(\alpha) < s_2(\alpha))$.

\bigskip 

Then ${\Bbb D}_F$ is a partial order and is $\kappa$-closed and
$\kappa^+$-cc (assuming $\kappa^{<\kappa} = \kappa$).
Let $G$ be ${\Bbb D}_F$-generic over $M$ and
$g = \cup \{ s : \exists x ((s,x) \in G)\}$.
Then since $G$ intersects suitably chosen dense subsets of 
${\Bbb D}_F$ in $M$, $g$ is a function from $\kappa$ to $\kappa$
which eventually dominates every function in $F$, i.e.
$\forall f \in F (f <^* g)$.

\bigskip

\noindent {\bf Definition 13.}\ \ Assume that ${\Bbb P} \in M$,
$( \hbox{ ${\Bbb P}$ is $\kappa$-closed} )^M$, and 
${\bf 1} \Vdash \hbox{ ``\ $\varphi \subseteq \kappa^\kappa$\ ''}$.
The standard ${\Bbb P}$-name for ${\Bbb D}_\varphi$ is
$\langle \psi, \leq_\psi, {\bf 1}_\psi \rangle$, where
$$\eqalign{\psi = \{
            \langle op(\check s, \phi ), {\bf 1}_{\Bbb P} \rangle :
	    &
  \hbox{
         $s$ is a function $\wedge$ ${\rm dom}(s) \in \kappa$ $\wedge$
          ${\rm ran}(s) \subseteq \kappa$ $\wedge$}\cr
          &\hbox{    ${\bf 1}_{\Bbb P} \Vdash
                    \hbox{
                           ``\ $\phi \subseteq \varphi \wedge
				 \mid \phi \mid < \kappa$\ ''
                                            }$ $\wedge$
            $\phi$ is a nice name for a subset of $\varphi$
        }
          \}\cr}$$
and ${\bf 1}_\psi = op(\check 0, \check 0)$.

\bigskip

The choice of the ${\Bbb P}$-name $\leq_\psi$ is, once again, irrelevant
as long as it is forced by ${\bf 1}_{\Bbb P}$ to be the correct 
partial order on ${\Bbb D}_\varphi$.

In $M$, let ${\Bbb P}$, $\varphi$, $\psi$, be as above.
Let $G$ be ${\Bbb P}$-generic over $M$ and $F = \varphi_G$.
Then, in $M[G]$, $\psi_G = {\Bbb D}_F$. In addition, $\psi$ is
full for $<\!\!\kappa$-sequences.
The use of full names for $<\!\!\kappa$-sequences will guarantee,
as indicated earlier, that the iteration is $\kappa$-closed.
The use of standard names will imply that the iteration also satisfies
the $\kappa^+$-cc so that all the cardinals are preserved in the 
final model.

Now follows the main result of this section. 

\bigskip

\noindent {\bf Theorem 14.}\ \ Let $N$ be a c.t.m. for $ZFC + GCH$
and, in $N$, let
$\kappa < \lambda \leq \mu \leq \theta$ be cardinals such that
$\kappa, \lambda, \mu$ are regular and ${\rm cf}(\theta) > \kappa$.
Then there is a cardinal preserving extension $M[G]$ of $N$
such that 
$$M[G] \models \hbox{`` ${\frak t}_\kappa = \lambda \wedge
			 {\frak b}_\kappa = \mu \wedge
			 2^\kappa = \theta$ ''}.$$

\bigskip

Proof:\ \ Let $\alpha$, $\beta$
be cardinals with $\alpha$ regular, $\alpha < \beta$, and
${\rm cf}(\beta) > \alpha$. Then 
${\Bbb Fn}(\beta \times \alpha, 2, \alpha)$ is the standard partial
order for adding $\beta$-many subsets of $\alpha$ (see [2]).
It is $\alpha$-closed and $\alpha^+$-cc 
(assuming $\alpha^{<\alpha} = \alpha$), so it preserves cardinals.

Let $N$ be a c.t.m. for $ZFC + GCH$. In $N$, let
$\kappa < \lambda \leq \mu \leq \theta$ be cardinals such that
$\kappa$, $\lambda$, $\mu$ are regular and ${\rm cf}(\theta) > \kappa$.
The goal is to produce an extension of $N$ in which
${\frak t}_\kappa = \lambda$, 
${\frak b}_\kappa = \mu$ and $2^\kappa = \theta$. Let $H$ be
${\Bbb Fn}( \theta^+\times \lambda, 2, \lambda)$-generic over $N$
and let $N[H] = M$. Then
$$M \models \hbox{``\ $ZFC + \forall \xi < \lambda (2^\xi = \xi^+) + 
   2^\lambda = \theta^+$\ ''}$$
$\kappa$, $\lambda$, $\mu$ are still regular and all the cardinals are 
preserved. Now, in $M$, perform an iterated forcing construction
of length $\theta \cdot \mu$ (ordinal product) with
$<\!\!\kappa$-supports, i.e. build an iterated forcing construction
$$\langle\langle {\Bbb P}_\xi : \xi\leq \theta \cdot \mu \rangle,
	 \langle \pi_\xi : \xi < \theta \cdot \mu \rangle\rangle$$
with supports of size less than $\kappa$.

Given ${\Bbb P}_\xi$, if $\xi$ is not of the form $\theta \cdot \xi$,
list all the ${\Bbb P}_\xi$-names for towers in
$({\cal P}(\kappa), \subset^*)$ of size $\eta$ for all
$\kappa < \eta < \lambda$; for example, let
$\langle \sigma^\xi_\gamma : \gamma < \theta \rangle$
enumerate all ${\Bbb P}_\xi$-names $\sigma$ such that for some $\eta$,
with $\kappa < \eta < \lambda$, $\sigma$ is a nice ${\Bbb P}_\xi$-name
for a subset of $(\eta \times \kappa)\check{\vphantom{o}}$
with the property that there is a name $\tau^\xi_\gamma$ such that
$${\bf 1}\Vdash \hbox{``\ $\tau^\xi_\gamma = \{ x \subseteq \kappa :
	\exists \zeta < \eta ( x = \{ \nu : (\zeta, \nu) \in 
	\sigma^\xi_\gamma\} ) \}$ is a tower in
	$({\cal P}(\kappa), \subset^*)$ of size $\eta$\ ''}.$$

Let 
$\Theta = (\theta \cdot \mu) \setminus \{ \theta \cdot \xi :
	       \xi < \mu \}$
and let 
$f : \Theta \to (\theta \cdot \mu) \times \theta$
be a bookkeeping function such that $f$ is onto and
$\forall \xi, \beta, \gamma (f(\xi) = (\beta, \gamma) \to
  \beta < \xi )$.
If $f(\xi) = (\beta, \gamma)$, let $\tau_\xi$  be a ${\Bbb P}_\xi$-name
for the same object for which $\tau^\beta_\gamma$ is a ${\Bbb P}_\beta$-name.
Let $\pi_\xi$ be the standard ${\Bbb P}_\xi$-name for
${\Bbb T}_{\tau_\xi}$. And if $\xi$ is of the form $\theta \cdot \zeta$,
let $\varphi_\xi$ be a ${\Bbb P}_\xi$-name
for $\kappa^\kappa$ and let $\pi_\xi$ be the standard ${\Bbb P}_\xi$-name
for ${\Bbb D}_{\varphi_\xi}$. This finishes the iteration.

By Lemma 9 ${\Bbb P}_{\theta \cdot \mu}$ is $\kappa$-closed in $M$.
In fact, ${\Bbb P}_{\theta \cdot \mu}$ has the property that each decreasing 
sequence of length $<\kappa$ has a greatest lower bound so that the set
${\Bbb P}'$ of elements $p\in {\Bbb P}_{\theta \cdot \mu}$ with the property
that the first coordinate of $p(\gamma)$, for $\gamma \in {\rm dom}(p)$,
is a real object and not just a ${\Bbb P}_\gamma$-name, is dense in
${\Bbb P}_{\theta \cdot \mu}$. Therefore, to show that
${\Bbb P}_{\theta \cdot \mu}$  also has the $\kappa^+$-cc in $M$
it suffices to show that ${\Bbb P}'$ has the $\kappa^+$-cc in $M$.
So, in $M$, let $p^\gamma \in {\Bbb P}'$
for $\gamma < \kappa^+$. By $\kappa^{<\kappa} = \kappa$, the 
$\Delta$-system lemma (see Theorem II 1.6 in [2])  implies that 
there is an $X \in [\kappa^+]^{\kappa^+}$ such that
$\{ {\rm support}(p^\gamma) : \gamma < \kappa^+ \}$
for a $\Delta$-system with root $r$. 
Let $p^\gamma = \langle \rho^\gamma_\xi : \xi < \theta \cdot \mu \rangle$,
and let $ \rho^\gamma_\xi = op ( \check s^\gamma_\xi, \sigma^\gamma_\xi )$.
By $\kappa^{<\kappa} = \kappa$, there is a
$Y \in [X]^{\kappa^+}$ such that for all $\xi \in r$, the $s^\gamma_\xi$
for $\gamma \in Y$ are all the same; say $s^\gamma_\xi = s_\xi$
for $\xi \in r$ and $\gamma \in Y$.
But then the $ p^\gamma$ for $\gamma \in Y$ are pairwise compatible;
to see this observe that if $\gamma, \delta \in Y$, then
$p^\gamma$, $p^\delta$ have as a common extension
$\langle p_\xi : \xi < \theta \cdot \mu \rangle$, where $\rho_\xi$ is

$(a)$\ \ $\rho^\gamma_\xi$ if $\xi \not \in {\rm support}(p^\delta)$,

$(b)$\ \ $\rho^\delta_\xi$ if $\xi \not \in {\rm support}(p^\gamma)$,

$(c)$\ \ $op(\check s_\xi, \sigma_\xi)$ if $\xi \in r$,

\noindent where $\sigma_\xi$ is a nice name which satisfies
${\bf 1}_\xi \Vdash \hbox{``\  $\sigma_\xi = 
          \sigma^\gamma_\xi \cup \sigma^\delta_\xi$\ ''}$.
So  ${\Bbb P}_{\theta \cdot \mu}$ has the $\kappa^+$-cc and
together with being $\kappa$-closed preserves all the cardinal
numbers. Let $G$ be ${\Bbb P}_{\theta \cdot \mu}$-generic
over $M$. Since at each stage of the form $\theta \cdot \xi$,
a function from $\kappa$ to $\kappa$ is added which eventually dominates
all the functions in $\kappa^\kappa$ constructed by that stage,
it follows that, in $M[G]$, there is a scale in $(\kappa^\kappa, <^*)$
of order type $\mu$ so that ${\frak b}_\kappa = \mu$. In addition,
since at each stage of the iteration a new element to
$\kappa^\kappa$ or ${\cal P}(\kappa)$ is added,
it follows that 
$M[G] \models \hbox{``\ $2^\kappa = 
		\mid \theta \cdot \mu \mid = \theta$\ ''}$.
Finally, $M[G]$ contains no towers in 
$({\cal P}(\kappa), \subseteq^*)$ of order type $\eta$ for
$\kappa < \eta < \lambda$ since by the bookkeeping device all such
towers are considered and eventually destroyed at some stage of 
the iteration, so that ${\frak t}_\kappa \geq \lambda$.
However,
$M[G] \models \hbox{``\ $\forall \xi ( \kappa \leq \xi < \lambda \to
    2^\xi = \theta)$\ ''}$ and
$M[G] \models \hbox{ ``\ $2^\lambda = \theta^+$\ ''}$ since
$M \models \hbox{``\ $2^\lambda = \theta^+$\ ''}$ so that
by the previous section ${\frak t}_\kappa = \lambda$.
This finishes the proof of  this theorem. \qquad $\square$

\bigskip

\noindent {\bf  References}

\bigskip

\noindent [1] J. Cummings and S. Shelah, Cardinal Invariants above the
              Continuum, Annals of Pure and Applied Logic, 75 no. 3 (1995)
              251--268.

\medskip

\noindent [2] K. Kunen, Set Theory. An introduction to independence proofs,
              North Holland (1980). 

\medskip

\noindent [3] F. Rothberger, Sur un ensemble toujours de premiere categorie
              qui est depourvu de la propri\' et\' e $\lambda$, Fundamenta 
              Mathematicae, 32 (1939) 294--300.

\medskip

\noindent [4] F. Rothberger, On some problems of Hausdorff and
              Sierpi\' nski, Fundamenta Mathematicae, 35 (1948) 29--46.

\medskip

\noindent [5] E. K. van Douwen, Integers in Topology, Handbook of 
              Set-Theoretic Topology, eds. K. Kunen and J.E. Vaughan,
              North Holland (1984).

\end